# A Small Sample Correction for Estimating Attributable Risk in Case-Control Studies

Daniel B. Rubin

University of California, Berkeley

October 29, 2018

The attributable risk, often called the population attributable risk, is in many epidemiological contexts a more relevant measure of exposure-disease association than the excess risk, relative risk, or odds ratio. When estimating attributable risk with case-control data and a rare disease, we present a simple correction to the standard approach making it essentially unbiased, and also less noisy. As with analogous corrections given in Jewell (1986) for other measures of association, the adjustment often won't make a substantial difference unless the sample size is very small or point estimates are desired within fine strata, but we discuss the possible utility for applications.

Key phrases: attributable risk, case-control studies, small sample corrections.

# 1 Attributable risk

Although an exposure may be strongly associated with a disease, an intervention removing the risk factor can produce a limited public health benefit if the exposure is uncommon. Herein lies the interpretational problem with the measure of exposure-disease association

$$RR = \text{relative risk} = \text{pr(disease|exposed)}/\text{pr(disease|unexposed)}.$$

A different measure of association is Levin's (1953) attributable risk

$$\begin{aligned} AR &= \left\{\text{pr(disease)} - \text{pr(disease|unexposed)}\right\}/\text{pr(disease)} \\ &= \text{pr(exposed|disease)}\left(1 - 1/RR\right), \end{aligned}$$

where equality of the two lines can be seen from (5.1) of Jewell (2003). The attributable risk is informally the proportion of cases that could be saved by eliminating exposure, assuming no confounding. For instance, a value of 0.1 suggests a tenth of cases are due to exposure. When contemplating interventions, such knowledge can be more relevant than the relative risk.

While the attributable risk is not directly identifiable with case-control sampling, for a rare disease it is approximately

$$AR^\star = \text{pr(exposed|disease)}\left(1 - 1/OR\right),$$

for

$$OR = \text{odds ratio} = \frac{\text{pr(exposed|disease)}}{\text{pr(unexposed|disease)}} \frac{\text{pr(unexposed|disease free)}}{\text{pr(exposed|disease free)}}.$$

The reasoning behind the approximation is that

$$AR^\star - AR = \frac{\text{pr(exposed|disease)}}{RR}\left(1 - \frac{\text{pr(disease free|exposed)}}{\text{pr(disease free|unexposed)}}\right),$$

and the factor on the right will be close to zero if the disease is uncommon for both the exposed and unexposed. Case-control studies are most appropriate for rare diseases, so the following discussion assumes $AR^\star \approx AR$, and we henceforth take $AR^\star$ as our desired measure of association.

In a case-control study, we assess exposure status for $m$ cases and $n$ controls. The statistical model is

$$\text{number of cases exposed} = A \sim \text{Bin}(m, q),\ 0 < q < 1$$
$$\text{number of controls exposed} = B \sim \text{Bin}(n, p),\ 0 < p < 1,\ A \perp B$$
$$\text{number of cases unexposed} = C = m - A$$
$$\text{number of controls unexposed} = D = n - B.$$

We can now formally define our parameter of interest as

$$AR^\star = q\left(1 - \frac{1-q}{q}\frac{p}{1-p}\right).$$

## 2 Small sample correction

For additional details we refer to Whittemore (1982) or Jewell (2003, section 7.4), but the usual maximum likelihood estimator is given by

$$\widehat{AR} = \frac{A}{A+C} - \frac{BC}{(A+C)D}.$$

There is ambiguity when $D = 0$, and if the definition is taken literally the estimator has bias of $-\infty$. Since $A+C = m$, this part of the two denominators does not cause problems. We propose the small sample correction

$$\widehat{AR}_{SS} = \frac{A}{A+C} - \frac{BC}{(A+C)(D+1)}.$$

The following theorem is proven in the appendix, and states that bias decreases exponentially with the number of controls.

**Theorem 1.** $E[\widehat{AR}_{SS}] - AR^\star = \left(\frac{1-q}{1-p}\right) p^{n+1} = O(e^{-n})$.

For problems in statistics where exact unbiasedness is impossible, exponentially decreasing bias is somewhat unusual. A Taylor expansion can frequently show ratios of empirical means and the like have bias $\Omega(n^{-1})$, where $n$ is the sample size. The jackknife was originally introduced to reduce such biases to $O(n^{-2})$ (Quenouille, 1956), which is a far cry from $O(e^{-n})$.

Ideally we would only investigate exposures that are associated with disease, or at least neutral. In such a setting, the $(1-q)/(1-p)$ constant in the bias is no larger than one.

In line with its interpretation of yielding the percentage of cases due to exposure, attributable risk is often reported to two decimal places. Hence, say our estimator is essentially unbiased if bias is less than 0.005. If the exposure is not ubiquitous, so less than half of controls are expected to be exposed, then only nine controls are needed to ensure essential unbiasedness.

Table 1 compares the bias of the standard and corrected estimator. As mentioned, the standard estimator does not even have finite bias, but we artificially level the playing field by setting $\widehat{AR} = \widehat{AR}_{SS}$ when $D = 0$. The exact bias of the standard estimator can then be calculated by enumerating outcomes for the binomial distribution. Results are shown for a handful of parameter values $(q, p)$ across the unit square, for a small study with ten cases and ten controls. The corrected estimator clearly seems preferable.

Table 1: Bias of the standard and corrected estimator, with $m = 10$ cases and $n = 10$ controls, where $\widehat{AR}$ is set to $\widehat{AR}_{SS}$ if $D = 0$. The corrected estimator appears to have smaller bias.

| $q$ | $p$ | $AR^\star$ | Bias($\widehat{AR}$) | Bias($\widehat{AR}_{SS}$) |
|---|---|---|---|---|
| 0.20 | 0.20 | 0.0000 | -0.0297 | 0.0000 |
|  | 0.40 | -0.3333 | -0.1216 | 0.0001 |
|  | 0.60 | -1.0000 | -0.4647 | 0.0073 |
|  | 0.80 | -3.0000 | -1.0640 | 0.3436 |
| 0.40 | 0.20 | 0.2500 | -0.0223 | 0.0000 |
|  | 0.40 | 0.0000 | -0.0912 | 0.0000 |
|  | 0.60 | -0.5000 | -0.3485 | 0.0054 |
|  | 0.80 | -2.0000 | -0.7980 | 0.2577 |
| 0.60 | 0.20 | 0.5000 | -0.0149 | 0.0000 |
|  | 0.40 | 0.3333 | -.0608 | 0.0000 |
|  | 0.60 | 0.0000 | -0.2324 | 0.0036 |
|  | 0.80 | -1.0000 | -0.5320 | 0.1718 |
| 0.80 | 0.20 | 0.7500 | -.0074 | 0.0000 |
|  | 0.40 | 0.6667 | -.0304 | 0.0000 |
|  | 0.60 | 0.5000 | -0.1162 | 0.0018 |
|  | 0.80 | 0.0000 | -0.2660 | 0.0860 |

The bias reduction would be a Pyrrhic victory if it induced a large spike in variance. Fortunately, this does not occur because the correction tends to stabilize things. The standard method has infinite variance since $D = 0$ can appear in the denominator, but even after conditioning on the estimator being properly defined, the small sample correction still improves variance. The following theorem makes this precise, and our proof is in the appendix.

**Theorem 2.** $var(\widehat{AR}_{SS}|D \neq 0) < var(\widehat{AR})|D \neq 0)$.

The correction is thus meant to reduce both bias and variance, the components of mean squared error, or the traditional measure of risk in statistical decision theory. Our problem is a counterexample to the common situation where bias correction is ill-advised since the variance increase is prohibitive (Doss and Sethuraman, 1989; Efron and Tibshirani, 1993, section 10.6).

## 3 Does the correction make a difference?

Our method has several theoretical advantages that we have discussed, it can be implemented by hand, and we avoid definitional problems when $D = 0$. The procedure could prove useful with a small sample size, or when confounding necessitates cross-tabulation and point estimates are desired within relatively fine strata. However, replacing $D$ by $D+1$ in the estimator will often result in a very minor perturbation, and in this sense our proposal brings to mind the longstanding but frequently meaningless debate over whether to multiply a sample variance by $n/(n-1)$.

To see the utility of our method, we should note that it is an analog of relative risk and odds ratio corrections of Jewell (1986). Adjustments follow from the fact that if $X \sim \text{Bin}(n,p)$, then $(n+1)/(X+1)$ and $X/(n-X+1)$ are corrected estimators of $1/p$ and $p/(1-p)$. Jewell (2003, section 7.1) states that the "principal value of a small sample adjustment is in its alerting us to situations where the sample size is small enough to have a noticeable impact on an estimator, thereby suggesting that large sample approximations may be suspicious." This seems like good advice for attributable risk estimation, for which our method should have a modest but positive impact.

# Appendix

*Proof of Theorem 1.* We begin by noting that since $B \sim \text{Bin}(n, p)$, we have

$$
\begin{aligned}
E\left[\frac{B}{D+1}\right] &= E\left[\frac{B}{n-B+1}\right] \\
&= \sum_{k=0}^{n} \frac{n!}{k!(n-k)!} p^k (1-p)^{n-k} \frac{k}{n-k+1} \\
&= \sum_{k=1}^{n} \frac{n!}{k!(n-k)!} p^k (1-p)^{n-k} \frac{k}{n-k+1} \\
&= \frac{p}{1-p} \sum_{k=1}^{n} \frac{n!}{(k-1)!(n-(k-1))!} p^{k-1} (1-p)^{n-(k-1)} \\
&= \frac{p}{1-p} \sum_{i=0}^{n-1} \frac{n!}{i!(n-i)!} p^i (1-p)^{n-i} \text{ for } i = k-1 \\
&= \frac{p}{1-p} (1 - \text{pr}(B=n)) \\
&= \frac{p}{1-p} - \frac{p}{1-p} p^n \\
&= \frac{p}{1-p} - \frac{p^{n+1}}{1-p}.
\end{aligned}
$$

Further, $E[A/(A+C)] = E[A/m] = q$ since $A \sim \text{Bin}(m, q)$. Likewise, $E[C/(A+B)] = 1-q$. By the independence of $(A, C)$ and $(B, D)$, we conclude that

$$
\begin{aligned}
E[\widehat{AR}_{SS}] &= E\left[\frac{A}{A+C}\right] - E\left[\frac{C}{A+C}\frac{B}{D+1}\right] \\
&= E\left[\frac{A}{A+C}\right] - E\left[\frac{C}{A+C}\right] E\left[\frac{B}{D+1}\right] \\
&= q - (1-q)\frac{p}{1-p} + \left(\frac{1-q}{1-p}\right) p^{n+1}.
\end{aligned}
$$

The desired result now follows after expressing

$$
AR^\star = q\left(1 - \frac{1-q}{q}\frac{p}{1-p}\right) = q - (1-q)\frac{p}{1-p}. \quad \Box
$$

*Proof of Theorem 2.* We first write

$$\widehat{AR} = \widehat{AR}_{SS} + \frac{A-m}{m}\left(\frac{n-D}{D} - \frac{n-D}{D+1}\right),$$

yielding

$$\begin{aligned}
\text{var}(\widehat{AR}|D \neq 0) &= \text{var}(\widehat{AR}_{SS}|D \neq 0) \\
&+ \text{var}\left[\frac{A-m}{m}\left(\frac{n-D}{D} - \frac{n-D}{D+1}\right)|D \neq 0\right] \\
&+ 2\,\text{cov}\left[\widehat{AR}_{SS}, \frac{A-m}{m}\left(\frac{n-D}{D} - \frac{n-D}{D+1}\right)|D \neq 0\right].
\end{aligned}$$

Both arguments of the covariance are increasing in $A$ and $D$. This implies the covariance term is positive, and hence the desired inequality. $\square$